\documentclass[a4paper,12pt]{article}

\usepackage[left=2cm,right=2cm, top=2cm,bottom=2.5cm,bindingoffset=0cm]{geometry}

\usepackage{verbatim}
\usepackage{amsmath}
\usepackage{amsthm}
\usepackage{amssymb}
\usepackage{delarray}
\usepackage{mathrsfs}
\usepackage{graphicx}
\usepackage{cite}
\usepackage{xcolor}

 \allowdisplaybreaks

\begin{document}

\begin{center}
{\large\bf INVERSE SPECTRAL PROBLEMS FOR HILL-TYPE OPERATORS WITH FROZEN ARGUMENT}
\\[0.5cm]

{\large\bf Sergey Buterin\footnote{Department of Mathematics, Saratov State University, Astrakhanskaya 83,
Saratov 410012, Russia, email: buterinsa@info.sgu.ru} and Yi-Teng~Hu\footnote{School of Mathematics and
Statistics, Xidian University, Xi'an, Shaanxi, 710071, China, email: ythu@xidian.edu.cn}} \\[0.2cm]
\end{center}

{\bf Abstract.} The paper deals with nonlocal differential operators possessing a term with frozen (fixed)
argument appearing, in particular, in modelling various physical systems with feedback. The presence of a
feedback means that the external affect on the system depends on its current state. If this state is taken
into account only at some fixed physical point, then mathematically this corresponds to an operator with
frozen argument. In the present paper, we consider the operator $Ly\equiv-y''(x)+q(x)y(a),$
$y^{(\nu)}(0)=\gamma y^{(\nu)}(1),$ $\nu=0,1,$ where $\gamma\in{\mathbb C}\setminus\{0\}.$ The operator $L$
is a nonlocal analog of the classical Hill operator describing various processes in cyclic or periodic media.
We study two inverse problems of recovering the complex-valued square-integrable potential $q(x)$ from some
spectral information about $L.$ The first problem involves only single spectrum as the input data. We obtain
complete characterization of the spectrum and prove that its specification determines $q(x)$ uniquely if and
only if $\gamma\ne\pm1.$ For the rest (periodic and antiperiodic) cases, we describe classes of iso-spectral
potentials and provide restrictions under which the uniqueness holds. The second inverse problem deals with
recovering $q(x)$ from the two spectra related to $\gamma=\pm1.$ We obtain necessary and sufficient
conditions for its solvability and establish that uniqueness holds if and only if $a=0,1.$ For $a\in(0,1),$
we describe classes of iso-bispectral potentials and give restrictions under which the uniqueness resumes.
Algorithms for solving both inverse problems are provided. In the appendix, we prove Riesz-basisness of an
auxiliary two-sided sequence of sines.


\\

{\large\bf 1. Introduction}
\\

In recent years, there appeared a considerable interest in nonlocal differential operators possessing a term
with frozen argument (see \cite{Nizh-09, AlbHryNizh, Nizh-10, Nizh-11, Nizh-12, BBV, BV, XY19-1, BK, HBY,
XY19-2, Hu20, Wang20} and references therein). Operators of this kind form an important and illustrative
class of the so-called loaded differential operators (see, e.g., \cite{Isk, Kral, Nakh76, NakhBor77, Nakh12,
Lom14, Lom15}), which often appear in mathematical physics. The corresponding loaded equations can be
characterized by the presence of a trace of the unknown function. Below we describe some typical models of
physical systems with feedback
leading to nonlocal differential operators with frozen argument. 
The presence of a feedback means that the external affect on the system depends on its current state. If this
state is taken into account only at some fixed physical point of the system, then mathematically this
corresponds to an operator with frozen argument.

We focus on the boundary value problem ${\cal L}(q(x),a,\gamma)$ of the form
\begin{equation}\label{1.1}
\ell y:=-y''(x)+q(x)y(a)=\lambda y(x), \quad 0<x<1,
\end{equation}
\begin{equation}\label{1.2}
y^{(\nu)}(0)=\gamma y^{(\nu)}(1), \quad \nu=0,1,
\end{equation}
where $q(x)\in L_2(0,1)$ is a complex-valued function, $\lambda$ is the spectral parameter and
$\gamma\in{\mathbb C}\setminus\{0\},$ while $a\in[0,1].$ The operator generated by the
functional-differential expression $\ell$ and equipped with boundary conditions (\ref{1.2}) is called {\it
Hill-type operator with frozen argument}. This operator is a nonlocal loaded analog of the classical Hill
operator appearing after separating variables in partial differential equations describing various processes
in cyclic or periodic media.

In the present paper, we study inverse spectral problems of recovering the potential $q(x).$ The most
complete results in the inverse spectral theory are known for purely differential operators (local), see
monographs \cite{Mar, Lev, FrYur-1
}. The inverse problems arise in mathematics, mechanics, physics, geophysics, electronics and other branches
of natural sciences and engineering, including nanoscale technology. In particular, for the classical Hill
operator inverse problems were studied in \cite{Mar,Stan,MO, ST, KK, Yur17, Lev} and other works. For
operators with frozen argument as well as for other types of nonlocal operators, classical methods of the
inverse spectral theory do not work.

Equation (\ref{1.1}) appears, for example, after applying the Fourier method of separation of variables to
the following loaded parabolic partial differential equation:
\begin{equation}\label{1.1.1}
\frac{\partial}{\partial t}u(x,t)=\frac{\partial^2}{\partial x^2}u(x,t)+f(x,t), \quad f(x,t)= -q(x)u(a,t),
\quad 0<x<1, \quad t>0,
\end{equation}
involving the trace $u(a,t)$ of the unknown function $u(x,t).$ Consider also the initial condition
\begin{equation}\label{1.1.2}
u(x,0)=\varphi(x), \quad 0<x<1,
\end{equation}
and the homogeneous boundary conditions
\begin{equation}\label{1.1.3}
\frac{\partial^\alpha}{\partial x^\alpha}u(x,t)\Big|_{x=0}=\frac{\partial^\beta}{\partial
x^\beta}u(x,t)\Big|_{x=1}=0, \quad t>0,
\end{equation}
where $\alpha,\beta\in\{0,1\}$ are fixed. It is well known that the initial-boundary value problem
(\ref{1.1.1})--(\ref{1.1.3}) models heat conduction in a rod of unit length possessing an external
distributed heat source $f(x,t).$ In our case, this heat source is described by the function $-q(x)u(a,t),$
i.e. its power is proportional to the temperature $u(a,t)$ at the fixed point~$a$ of the rod. Such a model
can be implemented by an electric conductive rod of a constant thermal conductivity but possessing the
variable electrical resistance $R(x)=-q(x)\ge0$ independent of the temperature of the rod. In~order to create
the heat source as in (\ref{1.1.1}), the voltage $U(t)$ being applied to the ends of the rod should be
proportional to $\sqrt{u(a,t)}:$
$$
U(t)={\cal R}\sqrt{u(a,t)}, \quad {\cal R}=-\int_0^1q(\xi)\,d\xi,
$$
where ${\cal R}$ is the full resistance of the rod. Indeed, the current $I(t)$ in the rod does not depend
on~$x$ and, by Ohm's law, we have $I(t)=U(t)/{\cal R}=\sqrt{u(a,t)},$ while the density of the created heat
source in each point $x\in[0,1]$ of the rod equals to $I^2(t)R(x)=-q(x)u(a,t)=f(x,t).$

Similarly, the loaded hyperbolic equation
\begin{equation}\label{1.1.4}
\frac{\partial^2}{\partial t^2}u(x,t)=\frac{\partial^2}{\partial x^2}u(x,t) -q(x)u(a,t), \quad 0<x<1, \quad
t>0,
\end{equation}
under the initial conditions
\begin{equation}\label{1.1.5}
u(x,0)=\varphi(x), \quad \frac{\partial}{\partial t}u(x,t)\Big|_{t=0}=\psi(x), \quad 0<x<1,
\end{equation}
as well as the boundary conditions (\ref{1.1.3}), with $q(x),$ $\varphi(x)$ and $\psi(x)$ belonging to some
appropriate classes, models an oscillatory process under a damping external force that is proportional in
each point $x$ both to the function $q(x)$ and to the displacement $u(a,t)$ of the oscillating system in the
fixed point $a$ equipped with a displacement sensor. An example occurs if a vibrating wire is affected by a
magnetic field exerting a force per unit mass represented by $-q(x)u(a,t),$ i.e. depending on the lateral
displacement $u(a,t)$ at the point $a$ at the time~$t.$

After separating variables, both above-mentioned initial-boundary value problems yield one and the same
eigenvalue problem ${\cal B}:={\cal B}(q(x),a,\alpha,\beta)$ consisting of the ordinary
functional-differential equation with frozen argument (\ref{1.1}) and the separated boundary conditions
\begin{equation}\label{1.2-sep}
y^{(\alpha)}(0)=y^{(\beta)}(1)=0.
\end{equation}
We note that some analogous as well as different models leading to ordinary functional-differential equations
with one or several frozen arguments were given, e.g., in \cite{Kral, Nakh12}.

Various aspects of inverse problems for operators with frozen argument were studied in \cite{Nizh-09,
AlbHryNizh, Nizh-10, Nizh-11, Nizh-12, BBV, BV, BK, HBY, XY19-2, Wang20}. In particular, in \cite{BK} the
problem ${\cal B}$ was considered for $a\in[0,1]\cap{\mathbb Q},$ i.e. $a=j/k$ with some mutually prime
integers $j$ and~$k.$ It was established that unique recoverability of the potential $q(x)$ from the spectrum
of ${\cal B}$ depends on the values $\alpha,\,\beta$ as well as on the parity~of~$k.$ This implies two cases:
{\it non-degenerate} and {\it degenerate} ones (see also \cite{BBV,BV}), depending on whether the inverse
problem is uniquely solvable or not, respectively. In the degenerate case, asymptotically $k$-th part of the
spectrum degenerates in the sense that each $k$-th eigenvalue carries no information on the potential. In
\cite{Wang20}, it was shown that each irrational $a\in(0,1)$ implies the non-degenerate case for any pair of
$\alpha,\beta\in\{0,1\}.$

Let $\{\lambda_n\}_{n\ge0}$ be the spectrum of the boundary value problem ${\cal L}(q(x),a,\gamma).$ We start
with the following inverse problem.

\medskip
{\bf Inverse Problem 1.} Given $\{\lambda_n\}_{n\ge0},$ $a$ and $\gamma;$ find $q(x).$

\medskip
For example, the case $\gamma=1$ corresponds to the initial-boundary value problem for the heat
equation~(\ref{1.1.1}) under the initial condition (\ref{1.1.2}) as well as the periodic boundary conditions
$$
u(0,t)=u(1,t), \quad \frac{\partial}{\partial x}u(x,t)\Big|_{x=0}=\frac{\partial}{\partial
x}u(x,t)\Big|_{x=1}, \quad t>0,
$$
which models heat conduction in a thin closed ring of unit length parameterized by the variable $x\in[0,1],$
whose values~$0$ and~$1$ correspond to one and the same physical point of the ring. Analogously, by using
equation (\ref{1.1.4}) along with the initial conditions (\ref{1.1.5}), one can model the corresponding
oscillatory process. Each of these two models hints that it is sufficient to study the case $a=0.$ Moreover,
this sufficiency remains also for any $\gamma\ne0$ (see Lemma~4 in Section~3). This property allows making no
distinguishing between rational and irrational cases, which gives a new quality comparing with the case of
separated boundary conditions (\ref{1.2-sep}).

For any $a\in[0,1]$ and $\gamma\ne0,$ we obtain complete characterization of the spectrum and prove that its
specification determines $q(x)$ uniquely if and only if $\gamma\ne\pm1.$ For $\gamma=\pm1,$ we establish that
precisely half of the spectrum degenerates, and describe classes of iso-spectral potentials. Moreover, we
provide restrictions on the potential under which the uniqueness resumes. The proof is constructive and gives
algorithms for solving Inverse Problem~1.

Further, we study recovering the potential $q(x)$ from the periodic and the antiperiodic spectra. For
$\alpha=0,1,$ we denote by $\{\lambda_{n,\alpha}\}_{n\ge0}$ the spectrum of the boundary value problem ${\cal
L}_\alpha(q(x),a):={\cal L}(q(x),a,(-1)^\alpha)$ and consider the following inverse problem.

\medskip
{\bf Inverse Problem 2.} Given $\{\lambda_{n,\alpha}\}_{n\ge0},$ $\alpha=0,1,$ and $a;$ find $q(x).$

\medskip
We prove a uniqueness theorem and obtain a constructive procedure for solving this inverse problem along with
necessary and sufficient conditions for its solvability. In particular, it is established that the solution
is unique if and only if $a\in\{0,1\}.$ For $a\in(0,1),$ we describe classes of iso-bispectral potentials and
provide restrictions under which the uniqueness holds. In the latter case, characterization of the spectra,
besides asymptotics, includes a restriction on the type of the sum of the characteristic functions being an
entire function of order $1/2.$

The paper is organized as follows. In the next section, we study the characteristic function of the problem
${\cal L}(q(x),a,\gamma)$ and derive the so-called main equation of the inverse problem. Therein, we also
obtain asymptotics of the spectrum. In Section~3, we study Inverse Problem~1, while Inverse Problem~2 is
investigated in Section~4. In Appendix~A, we prove Riesz-basisness of one auxiliary functional system.
\\

{\large\bf 2. Characteristic function. Main equation}
\\

Consider the solutions $C(x,\lambda)$ and $S(x,\lambda)$ of equation (\ref{1.1}) under the initial conditions
$$
C(a,\lambda)=S'(a,\lambda)=1, \quad S(a,\lambda)=C'(a,\lambda)=0.
$$
Having put $\rho^2=\lambda,$ we arrive at the representations
\begin{equation}\label{2.1}
C(x,\lambda)=\cos\rho(x-a)+\int_a^x\frac{\sin\rho(x-t)}{\rho}q(t)\,\mathrm{d}t, \quad S(x,\lambda)=\frac{\sin\rho(x-a)}{\rho}.
\end{equation}
Clearly, eigenvalues of the problem (\ref{1.1}), (\ref{1.2}) coincide with zeros of the entire function
\begin{equation}\label{2.3}
\Delta(\lambda)=
\begin{vmatrix}
C(0,\lambda)-\gamma C(1,\lambda) & S(0,\lambda)-\gamma S(1,\lambda) \\
C'(0,\lambda)-\gamma C'(1,\lambda) & S'(0,\lambda)-\gamma S'(1,\lambda)
\end{vmatrix},
\end{equation}
which is called {\it characteristics function}. Consider the Wronski-type determinant $W(x,\lambda):=\langle
C(x,\lambda),S(x,\lambda)\rangle,$ where $\langle y(x),z(x)\rangle=y(x)z'(x)-y'(x)z(x).$ Using (\ref{2.1}), it is easy to calculate
\begin{equation} \label{W}
W(x,\lambda)=1-\int_0^{a-x} q(a-t)\frac{\sin\rho t}\rho\,dt.
\end{equation}
Unlike the Wronskian for the classical Sturm--Liouville equation, the function $W(x,\lambda)$ depends on $x$
and, moreover, may vanish for some values of $x\in[0,1].$ However, we need the designation $W(x,\lambda)$
only for brevity. The next lemma gives a fundamental representation for $\Delta(\lambda).$

\medskip
{\bf Remark 1.} Since the spectra of the problems ${\cal L}(q(x),a,\gamma)$ and ${\cal
L}(q(1-x),1-a,\gamma^{-1}),$ obviously, coincide, without loss of generality, one can assume that
$a\in[0,1/2].$

\medskip
{\bf Lemma 1. }{\it For any $\gamma,$ the characteristic function $\Delta(\lambda)$ of the problem
(\ref{1.1}), (\ref{1.2}) has the form
\begin{equation} \label{FundRep}
\Delta(\lambda)=1+\gamma^2 -2\gamma\cos\rho -\int_0^1 w(x)\frac{\sin\rho x}\rho\,dx,\quad w(x)\in L_2(0,1).
\end{equation}
Moreover, the following representation holds:
\begin{equation}\label{w}
w(x)=\left\{
\begin{array}{l}
\displaystyle \gamma^2 q(a+x) +q(a-x), \quad x\in(0,a),\\[3mm]
\displaystyle \gamma q(a+1-x) +\gamma^2 q(a+x), \quad x\in(a,1-a),\\[3mm]
\displaystyle \gamma(q(a+1-x)+q(a-1+x)), \quad x\in(1-a,1),\end{array}\right.
\end{equation}
where, without loss of generality, we assumed $2a\le1.$}

\medskip
{\it Proof.} According to (\ref{2.3}) and the definition of $W(x,\lambda),$ we calculate
\begin{equation}\label{DeltaW}
\Delta(\lambda)=W(0,\lambda) +\gamma\Delta_{1,0}(\lambda) -\gamma\Delta_{0,1}(\lambda) +\gamma^2
W(1,\lambda),
\end{equation}
where
$$
\Delta_{\alpha,\beta}(\lambda)=
\begin{vmatrix}
C^{(\alpha)}(0,\lambda) & S^{(\alpha)}(0,\lambda) \\
C^{(\beta)}(1,\lambda) & S^{(\beta)}(1,\lambda)
\end{vmatrix}.
$$
By virtue of Lemma~1 in \cite{BK}, for $\alpha\ne\beta$ we have
\begin{equation}\label{AlNeBe}
\Delta_{\alpha,\beta}(\lambda)=(-1)^\alpha\cos\rho+\int_0^1 w_{\alpha,\beta}(x)\frac{\sin\rho x}{\rho}\,dx,
\quad w_{\alpha,\beta}(x)\in L_2(0,1),
\end{equation}
where (for $2a\le1)$
\begin{equation}\label{wAB}
\displaystyle w_{\alpha,\beta}(x)=\frac12 \left\{
\begin{array}{cl}
\displaystyle  q(1-a+x) -q(1-a-x), & x\in(0,a),\\[3mm]
\displaystyle (-1)^{1+\beta}q(1+a-x) -q(1-a-x), & x\in(a,1-a),\\[3mm]
\displaystyle (-1)^{1+\beta}\Big(q(1+a-x)+q(x-1+a)\Big), & x\in(1-a,1).
\end{array}\right.
\end{equation}
Using (\ref{W}) and (\ref{DeltaW})--(\ref{wAB}), we arrive at (\ref{FundRep}) and (\ref{w}). $\hfill\Box$

\medskip
After assuming $w(x)$ to be known, relation (\ref{w}) can be considered as a linear functional equation with
respect to $q(x),$ which we refer to as {\it main equation} of the inverse problem.

For $\gamma\ne0,$ denote
\begin{equation}\label{alpha}
\alpha:=\frac1\pi\arccos\frac{1+\gamma^2}{2\gamma}\in\{z:{\rm
Re}\,z\in[0,1], \,{\rm Im}\,z\ge0\} \cup\{z:{\rm Re}\,z\in(0,1), \,{\rm Im}\,z<0\}.
\end{equation}
The following theorem describes behavior of the spectrum.

\medskip
{\bf Theorem 1.} {\it The spectrum $\{\lambda_n\}_{n\ge0}$ of the problem ${\cal L}(q(x),a,\gamma)$ has the
form
\begin{equation}\label{a=0}
\lambda_{2k}=(2 k+\alpha)^2\pi^2 +\varkappa_{2k}, \;\; k\ge0, \quad \lambda_{2k-1}=(2
k-\alpha)^2\pi^2+\varkappa_{2k-1}, \;\; k\ge1, \quad \{\varkappa_n\}_{n\ge0}\in l_2.
\end{equation}
Moreover, if $\gamma=\pm1$ (i.e. if $\alpha=0,1$), then
\begin{equation}\label{a=1}
\varkappa_{2k-1}=0, \quad k\ge1,
\end{equation}
i.e. in both the periodic and the antiperiodic cases, the spectrum degenerates in the sense of~(\ref{a=1}).}

\medskip
{\it Proof.} Rewrite (\ref{FundRep}) in the form
\begin{equation} \label{FundRep_2}
\Delta(\lambda)=4\gamma\sin\frac{\rho+\pi\alpha}2\sin\frac{\rho-\pi\alpha}2 -\int_0^1 w(x)\frac{\sin\rho
x}\rho\,dx, \quad w(x)\in L_2(0,1).
\end{equation}
Thus, by the known method involving Rouch\'e's theorem (see, e.g., \cite{FrYur-1}) one can prove that any
function of the form (\ref{FundRep}) has infinitely many zeros $\lambda_n,$ $n\ge0,$ which with account of
multiplicities have the form
\begin{equation}\label{eps}
\lambda_n=\rho_n^2, \quad\rho_{2k}=(2 k+\alpha)\pi +\varepsilon_{2k}, \;\; k\ge0, \;\; \rho_{2k-1}=(2
k-\alpha)\pi +\varepsilon_{2k-1}, \;\; k\ge1, \quad \varepsilon_n=o(1).
\end{equation}
Further we consider two cases.

(i) Let $\gamma\ne\pm1,$ i.e. $\sin\pi\alpha\ne0.$ Substituting (\ref{eps}) into (\ref{FundRep_2}), we arrive
at the relation
$$
4\gamma\Big((-1)^n\sin\pi\alpha\cos\frac{\varepsilon_n}2 +\frac{1+\gamma^2}{2\gamma}
\sin\frac{\varepsilon_n}2\Big)\sin\frac{\varepsilon_n}2 =\int_0^1 w(x)\frac{\sin\rho_n x}{\rho_n}\,dx.
$$
Since $\sin\rho=\rho+O(\rho^3)$ and $\cos\rho=1+O(\rho^2)$ as $\rho\to0,$ we refine
$\{n\varepsilon_n\}_{n\ge0}\in l_2,$ which along with (\ref{eps}) gives (\ref{a=0}) for $\gamma\ne\pm1.$

(ii) Let $\gamma\in\{-1,1\}.$ Then, by virtue of (\ref{w}), we have
\begin{equation}\label{par}
w(x)=(-1)^\frac{1-\gamma}2 w(1-x).
\end{equation}
For $\gamma=1,$ this implies
$$
\int_0^1 w(x)\sin\rho x\,dx =\int_0^\frac12 w(x)\Big(\sin\rho x +\sin\rho(1-x)\Big)\,dx =
2\sin\frac\rho2\int_0^\frac12 w\Big(\frac12-x\Big)\cos\rho x\,dx.
$$
Since $\alpha=0,$ we rewrite (\ref{FundRep_2}) in the form
\begin{equation}\label{par0}
\Delta(\lambda) =\frac2\rho\sin\frac\rho2 \Big(2\rho\sin\frac\rho2 -\int_0^\frac12
w\Big(\frac12-x\Big)\cos\rho x\,dx\Big),
\end{equation}
which yields (\ref{a=0}) and (\ref{a=1}) for $\gamma=1.$ Similarly, for $\gamma=-1,$ formula (\ref{par})
gives
$$
\int_0^1 w(x)\sin\rho x\,dx =\int_0^\frac12 w(x)\Big(\sin\rho x -\sin\rho(1-x)\Big)\,dx =
-2\cos\frac\rho2\int_0^\frac12 w\Big(\frac12-x\Big)\sin\rho x\,dx,
$$
which along with (\ref{FundRep_2}) and $\alpha=1$ implies
$$
\Delta(\lambda) =2\cos\frac\rho2 \Big(2\cos\frac\rho2 +\int_0^\frac12 w\Big(\frac12-x\Big)\frac{\sin\rho
x}\rho\,dx\Big).
$$
The latter formula gives (\ref{a=0}) and (\ref{a=1}) for $\gamma=-1.$ $\hfill\Box$

\medskip
The characteristic function is uniquely determined by the spectrum. In more detail, the following lemma hods.

\medskip
{\bf Lemma 2. }{\it Any function $\Delta(\lambda)$ of the form (\ref{FundRep}) is uniquely determined by its
zeros $\{\lambda_n\}_{n\ge0}$. Moreover, the following representation holds:
\begin{equation}\label{prod0}
\Delta(\lambda) =\left\{\begin{array}{cc} \displaystyle (1-\gamma)^2\prod_{n=0}^\infty
\frac{\lambda_n-\lambda}{\lambda_{n,\alpha}^0}, & \gamma\ne1,\\[3mm]
\displaystyle (\lambda-\lambda_0)\prod_{n=1}^\infty\frac{\lambda_n-\lambda}{\lambda_{n,\alpha}^0}, &
\gamma=1,
\end{array}\right.
\end{equation}
where
\begin{equation}\label{prod0-1}
\lambda_{n,\alpha}^0=(\rho_{n,\alpha}^0)^2, \;\; n\ge0, \quad \rho_{2k,\alpha}^0=(2 k+\alpha)\pi, \;\; k\ge0,
\quad \rho_{2k-1,\alpha}^0=(2 k-\alpha)\pi, \;\; k\ge1.
\end{equation}
}

{\it Proof.} By virtue of Hadamard's factorization theorem (see, e.g., \cite{BFY14}), formula (\ref{FundRep})
implies
\begin{equation}\label{4.14}
\Delta(\lambda)=C\lambda^s\prod_{\lambda_n\ne0}\Big(1-\frac\lambda{\lambda_n}\Big),
\end{equation}
where $C$ is some constant, while $s$ is the algebraic multiplicity of the null zero $\lambda_n=0.$ In
particular, we have
\begin{equation}\label{4.15}
\Delta_0(\lambda):=1+\gamma^2 -2\gamma\cos\rho=\left\{\begin{array}{rc} \displaystyle
(1-\gamma)^2\prod_{n=0}^\infty \Big(1-
\frac\lambda{\lambda_{n,\alpha}^0}\Big), & \gamma\ne1,\\[3mm]
\displaystyle \lambda\prod_{n=1}^\infty\Big(1- \frac\lambda{\lambda_{n,\alpha}^0}\Big), & \gamma=1.
\end{array}\right.
\end{equation}
Let $\gamma\ne1.$ Then, dividing (\ref{4.15}) by (\ref{4.14}), we obtain
\begin{equation}\label{4.16}
\frac{\Delta_0(\lambda)}{\Delta(\lambda)}=\frac{(1-\gamma)^2}C \prod_{\lambda_n=0}\Big(\frac1\lambda
-\frac1{\lambda_{n,\alpha}^0}\Big) \prod_{\lambda_n\ne0}\frac{\lambda_n}{\lambda_{n,\alpha}^0}
\prod_{\lambda_n\ne0}\frac{\lambda_{n,\alpha}^0-\lambda}{\lambda_n-\lambda},
\end{equation}
while (\ref{FundRep}) and (\ref{4.15}) imply $\Delta_0(\lambda)/\Delta(\lambda)\to1$ as $\lambda\to-\infty,$
which along with (\ref{4.16}) gives
$$
C=(1-\gamma)^2(-1)^s\prod_{\lambda_n=0}\frac1{\lambda_{n,\alpha}^0}
\prod_{\lambda_n\ne0}\frac{\lambda_n}{\lambda_{n,\alpha}^0}.
$$
Substituting this into (\ref{4.14}), we arrive at (\ref{prod0}) for $\gamma\ne1.$ The case $\gamma=1$ is
treated similarly. $\hfill\Box$

\medskip
For proving solvability of the inverse problems we will need also the following lemma.

\medskip
{\bf Lemma 3. }{\it Let $\gamma\ne\pm1.$ Then for any complex sequence $\{\lambda_n\}_{n\ge0}$ of the
form~(\ref{a=0}), the function $\Delta(\lambda)$ constructed by formula (\ref{prod0}) has the form
(\ref{FundRep}) with some function $w(x)\in L_2(0,1).$

Let $\gamma\in\{-1,1\}$ (i.e. $\alpha\in\{0,1\}$). Then for any complex sequence $\{\lambda_n\}_{n\ge0}$ of
the form~(\ref{a=0}) and satisfying the degeneration condition (\ref{a=1}), the function $\Delta(\lambda)$
constructed by formula~(\ref{prod0}) has the form~(\ref{FundRep}) with some function $w(x)\in L_2(0,1)$
obeying~(\ref{par}).}

\medskip
{\it Proof.} Let us first consider the case $\gamma=1,$ in which one should actually prove
representation~(\ref{par0}), where we have
$$
\frac2\rho\sin\frac\rho2 =\prod_{k=1}^\infty\Big(1-\frac\lambda{(2\pi k)^2}\Big).
$$
Thus, according to (\ref{a=0}) with $\alpha=0$ and (\ref{a=1}), it is sufficient to prove the relation
$$
(\lambda-\lambda_0)\prod_{k=1}^\infty\frac{\lambda_{2k}-\lambda}{(2\pi k)^2} =2\rho\sin\frac\rho2
-\int_0^\frac12 w\Big(\frac12-x\Big)\cos\rho x\,dx,
$$
which, in turn, can be established similarly to Lemma~3.3 in \cite{But07} or obtained as its corollary.

The case $\gamma=-1$ can be treated analogously. Assume now that $\gamma\ne\pm1$ and let us show that
$\{\rho_{n,\alpha}^0\Delta(\lambda_{n,\alpha}^0)\}_{n\ge0}\in l_2.$ Indeed, according to (\ref{prod0}),
(\ref{prod0-1}) and (\ref{4.15}), we have
$$
\Delta(\lambda)=\Delta_0(\lambda)\prod_{k=0}^\infty \frac{\lambda_k-\lambda}{\lambda_{k,\alpha}^0-\lambda}
=\frac{\lambda_n-\lambda}{\rho_{n,\alpha}^0+\rho} \cdot
\frac{\Delta_0(\lambda)}{{\rho_{n,\alpha}^0-\rho}}\prod_{{k\ne n}\atop{k=0}}^\infty
\frac{\lambda_k-\lambda}{\lambda_{k,\alpha}^0-\lambda}.
$$
Substituting $\lambda=\lambda_{n,\alpha}^0$ therein and using (\ref{a=0}), we get
$\rho_{n,\alpha}^0\Delta(\lambda_{n,\alpha}^0)= a_nb_n\varkappa_n,$ where
$$
a_n=\frac12\lim_{\rho\to\rho_{n,\alpha}^0}
\frac{\Delta_0(\lambda)}{\rho_{n,\alpha}^0-\rho}
=\gamma(-1)^{n+1}\sin\alpha\pi, \quad b_n=\prod_{{k\ne n}\atop{k=0}}^\infty
\frac{\lambda_k-\lambda_{n,\alpha}^0}{\lambda_{k,\alpha}^0-\lambda_{n,\alpha}^0} =\prod_{{k\ne
n}\atop{k=0}}^\infty \Big(1+\frac{\varkappa_k}{\lambda_{k,\alpha}^0-\lambda_{n,\alpha}^0}\Big).
$$
Using (\ref{alpha}), it is easy to show that $|\lambda_{k,\alpha}^0 -\lambda_{n,\alpha}^0|\ge C_\gamma(k+1)$
for $k\ne n,$ where $C_\gamma>0$ depends only on $\gamma.$ Hence, $|b_n|\le C.$ Further, note that the
functional system $\{\sin\rho_{n,\alpha}^0x\}_{n\ge0}$ coincides up to signs with the system
$\{\sin(2n+\alpha)\pi x\}_{n\in{\mathbb Z}},$ which, in turn, forms a Riesz basis in $L_2(0,1)$ as soon as
$\alpha\notin{\mathbb Z}.$ Indeed, for real non-integer $\alpha$'s this fact was proved in \cite{Bond18},
while for the general case the corresponding proof is provided in Appendix~A (see Lemma~A1). Thus, there
exists a unique function $w(x)\in L_2(0,1)$ such that
$$
\Delta(\lambda_{n,\alpha}^0)= -\int_0^1 w(x)\frac{\sin\rho_{n,\alpha}^0 x}{\rho_{n,\alpha}^0}\,dx,\quad
n\ge0.
$$
Consider the function
$$
\theta(\lambda)= -\int_0^1 w(x)\frac{\sin\rho x}{\rho}\,dx.
$$
It remains to show that $\Delta(\lambda)=\Delta_0(\lambda) +\theta(\lambda).$ For this purpose, we introduce
the function
$$
\sigma(\lambda)=\frac{\Delta(\lambda) -\Delta_0(\lambda)-\theta(\lambda)}{\Delta_0(\lambda)}
=F(\lambda)-1-\frac{\theta(\lambda)}{\Delta_0(\lambda)}, \quad F(\lambda)=\prod_{k=0}^\infty
\frac{\lambda_k-\lambda}{\lambda_{k,\alpha}^0-\lambda}.
$$
After removing singularities, the function $\sigma(\lambda)$ is entire in $\lambda.$ As in the proof of
Lemma~3.3 in \cite{But07}, one can show that $|F(\lambda)|<C_\delta$ for $\lambda\in
G_\delta:=\{\lambda=\rho^2:|\rho-(2n+\alpha)\pi|\ge \delta, n\in{\mathbb Z}\},$ $\delta>0,$ and
$F(\lambda)\to1$ as $\lambda\to-\infty.$ Moreover, we have $\rho\theta(\lambda)=o(\Delta_0(\lambda))$ in
$G_\delta$ as soon as $|\lambda|\to\infty.$ Thus, by virtue of the maximum modulus principle, the function
$\sigma(\lambda)$ is bounded and consequently, according to Liouville's theorem, it is constant. Since
$\sigma(\lambda)\to 0$ as $\lambda\to-\infty,$ we get $\sigma(\lambda)\equiv0,$ which finishes the proof.
$\hfill\Box$
\\

{\large\bf 3. Recovering from one spectrum}
\\

As mentioned in Introduction, when $\gamma=1,$ studying Inverse Problem~1 for any $a$ can, obviously, be
reduced to the case $a=0.$ The following lemma reveals this possibility also for any other $\gamma\ne0.$

\medskip
{\bf Lemma 4. }{\it For any $a\in(0,1],$ the spectrum of the problem ${\cal L}(q(x),a,\gamma)$ coincides with
the spectrum of ${\cal L}(q_a,0,\gamma),$ where the function $q_a(x)$ is determined by the formula
\begin{equation}\label{qa}
q_a(x)=\left\{
\begin{array}{l}
q(x+a), \quad x\in(0,1-a),\\[2mm]
\displaystyle\frac1\gamma q(x+a-1), \quad x\in(1-a,1).\end{array}\right.
\end{equation}}

\medskip
{\it Proof.} Denoting
\begin{equation}\label{ya}
y_a(x)=\left\{
\begin{array}{l}
y(x+a), \quad x\in[0,1-a],\\[2mm]
\displaystyle\frac1\gamma y(x+a-1), \quad x\in(1-a,1],\end{array}\right.
\end{equation}
we note that, by virtue of the boundary conditions (\ref{1.2}), $y_a(x)\in W_2^2[0,1].$ Moreover, since
$y(x)\in W_2^2[0,1],$ we have
\begin{equation}\label{1.2-ya}
y_a^{(\nu)}(0)=\gamma y_a^{(\nu)}(1), \quad \nu=0,1.
\end{equation}
Further, using (\ref{qa}) and (\ref{ya}), we rewrite equation (\ref{1.1}) in the form
\begin{equation}\label{1.1-ya}
-y_a''(x)+q_a(x)y_a(0)=\lambda y_a(x), \quad 0<x<1.
\end{equation}
It remains to note that one and the same $\lambda$ can be eigenvalue of the boundary value problem
(\ref{1.1}), (\ref{1.2}) and the problem (\ref{1.2-ya}), (\ref{1.1-ya}) only simultaneously and only of the
same multiplicity. $\hfill\Box$

\medskip
Now we proceed directly to studying Inverse Problem~1. Here and in the subsequent section, along with the
boundary value problem ${\cal L}:={\cal L}(q(x),a,\gamma)$ we consider a problem $\tilde{\cal L}:={\cal
L}(\tilde q(x),a,\gamma)$ of the same form but with a different potential $\tilde q(x).$ We agree that if a
certain symbol $\beta$ denotes an object related to the problem ${\cal L},$ then this symbol with tilde
$\tilde\beta$ will denote the corresponding object related to $\tilde {\cal L}.$ The following uniqueness
theorem holds.

\medskip
{\bf Theorem 2. }{\it Let $\gamma\ne\pm1.$ Then coincidence of the spectra: $\{\lambda_n\}_{n\ge0}
=\{\tilde\lambda_n\}_{n\ge0}$ implies coincidence of the potentials: $q(x)=\tilde q(x)$ a.e. on $(0,1).$

Let $\gamma\in\{-1,1\}.$ Assume that there exists an operator $K:L_2(0,1/2)\to L_2(0,1/2)$ with injective
$I+\gamma K,$ such that
\begin{equation}\label{dd:4-1}
q_a\Big(\frac12-x\Big)=K\Big(q_a\Big(\frac12+x\Big)\Big), \quad \tilde q_a\Big(\frac12-x\Big)=K\Big(\tilde
q_a\Big(\frac12+x\Big)\Big) \;\; {\rm a.e.\;\; on} \;\; \Big(0,\frac12\Big),
\end{equation}
where $I$ is the identity operator and the function $q_a(x)$ is determined by formula (\ref{qa}). Then
coincidence of the spectra: $\{\lambda_n\}_{n\ge0}=\{\tilde\lambda_n\}_{n\ge0}$ implies $q(x)=\tilde q(x)$
a.e. on $(0,1).$}

\medskip
{\it Proof.} Let  $\{\lambda_n\}_{n\ge0}=\{\tilde\lambda_n\}_{n\ge0}.$ Then, by virtue of (\ref{FundRep}) and
(\ref{prod0}), we always have $w(x)=\tilde w(x)$ a.e. on $(0,1).$ According to Lemmas~1 and~4, the main
equation (\ref{w}) takes the form
\begin{equation}\label{w2}
w(x)=\gamma q_a(1-x) +\gamma^2q_a(x), \quad x\in(0,1),
\end{equation}
where $q_a(x)$ is determined by (\ref{qa}). This can also be checked directly using (\ref{w}) and~(\ref{qa}).
Clearly, equation (\ref{w2}) is equivalent to the linear system
\begin{equation}\label{w3}
\left.\begin{array}{c} w(x)=\gamma q_a(1-x) +\gamma^2q_a(x),\\[3mm]
w(1-x)=\gamma^2 q_a(1-x) +\gamma q_a(x),
\end{array}\right\} \quad x\in\Big(0,\frac12\Big),
\end{equation}
whose determinant equals to $\gamma^2(1-\gamma^2).$ Thus, the system (\ref{w3}) is non-degenerate if and only
if $\gamma\ne\pm1.$ Hence, assuming $\gamma\ne\pm1$ and inverting (\ref{qa}), we get $q(x)=\tilde q(x)$ a.e.
on $(0,1).$

Let $\gamma\in\{-1,1\}.$ Then, with accordance to (\ref{par}), the equations in (\ref{w3}) are equivalent.
Making in the first equation of (\ref{w3}) the change of variable $x\to1/2-x$ and using (\ref{dd:4-1}), we
get
\begin{equation}\label{dd:4-3}
w\Big(\frac12-x\Big)=\gamma q_a\Big(\frac12+x\Big)+K\Big(q_a\Big(\frac12+x\Big)\Big), \;\;
w\Big(\frac12-x\Big)=\gamma\tilde q_a\Big(\frac12+x\Big)+K\Big(\tilde q_a\Big(\frac12+x\Big)\Big)
\end{equation}
a.e. on $(0,1/2).$ Since $I+\gamma K$ is injective, we get $q_a(x)=\tilde q_a(x)$ a.e. on $(1/2,1).$ Finally,
using (\ref{dd:4-1}) again, we arrive at $q_a(x)=\tilde q_a(x)$ a.e. on $(0,1/2),$ which finishes the proof.
$\hfill\Box$

\medskip
{\bf Remark 2.} Using (\ref{alpha}) and (\ref{a=0}), it is easy to show that specification of the spectrum
$\{\lambda_n\}_{n\ge0}$ determines also the constant~$\gamma$ uniquely if $\gamma\in\{-1,1\}$ and up to
inversion if $\gamma\ne\pm1.$

\medskip
{\bf Remark 3.} Clearly, the second part of Theorem~2 remains true also if one applies the operator $K$ to
the left-hand sides of the equalities in (\ref{dd:4-1}) instead of the right-hand ones.

For $\gamma\in\{1,-1\},$ condition \eqref{dd:4-1} may mean, in particular, the evenness (for $\gamma=1)$ or
the oddness (for $\gamma=-1)$ of the function $q_a(x)$ with respect to the midpoint of the interval $(0,1),$
i.e. $q_a(1/2-x)=\gamma q_a(1/2+x)$, $0<x<1/2.$ However, the case $K=-\gamma I$ is not covered by condition
\eqref{dd:4-1} and not eligible. Indeed, according to (\ref{dd:4-3}), in this case the spectrum of ${\cal
L}(q(x),a,\gamma)$ coincides with the one of ${\cal L}(0,a,\gamma)$ and, hence, carries no information on the
potential $q(x).$ The case of constant $K$ (i.e. when $K(f)$ is independent of $f)$ corresponds to a priori
specification of $q_a(x)$ on the subinterval $(0,1/2).$

\medskip
The next theorem means that Theorem~1 gives complete characterization of the spectrum.

\medskip
{\bf Theorem 3. }{\it Let $\gamma\ne\pm1.$ Then for any fixed $a\in[0,1]$ and any sequence of complex numbers
$\{\lambda_n\}_{n\ge0}$ of the form (\ref{a=0}), there exists a unique function $q(x)\in L_2(0,1)$ such that
$\{\lambda_n\}_{n\ge0}$ is the spectrum of the corresponding boundary value problem ${\cal
L}(q(x),a,\gamma).$

Let $\gamma\in\{-1,1\}.$ Then for any fixed $a\in[0,1]$ and any complex sequence $\{\lambda_n\}_{n\ge0}$ of
the form~(\ref{a=0}), obeying the degeneration condition (\ref{a=1}), there exists $q(x)\in L_2(0,1)$ (not
unique) such that  $\{\lambda_n\}_{n\ge0}$ is the spectrum of the corresponding problem ${\cal
L}(q(x),a,\gamma).$}

\medskip
{\it Proof.} Using $\{\lambda_n\}_{n\ge 0},$ we construct the function $\Delta(\lambda)$ by formula
(\ref{prod0}). By virtue of Lemma~3, it has the representation (\ref{FundRep}) with a certain function
$w(x)\in L_2(0,1).$ If $\gamma\in\{-1,1\},$ then this $w(x)$ additionally obeys condition (\ref{par}). Thus,
in any case, the linear system~(\ref{w3}) is consistent. Consider some its solution $q_a(x)\in L_2(0,1),$
which is not unique if and only if $\gamma\in\{-1,1\}.$  Find the function $q(x)\in L_2(0,1)$ determined by
formula (\ref{qa}) and consider the corresponding problem ${\cal L}:={\cal L}(q(x),a,\gamma).$ Let us show
that $\{\lambda_n\}_{n\ge 0}$ is its spectrum.

Indeed, by virtue of Lemma~4, the spectrum of ${\cal L}$ coincides with the one of the problem ${\cal
L}_a:={\cal L}(q_a,0,\gamma),$ where the function $q_a(x)$ is determined by formula (\ref{qa}). Let
$\Delta_a(\lambda)$ be the characteristic function of ${\cal L}_a.$ According to Lemma~1, it has the
representation
$$
\Delta_a(\lambda)=1+\gamma^2 -2\gamma\cos\rho -\int_0^1 w_a(x)\frac{\sin\rho x}\rho\,dx,
$$
where $w_a(x)=\gamma q_a(1-x) +\gamma^2q_a(x),$ $x\in(0,1).$ Comparing this with (\ref{w2}), we get
$w_a(x)=w(x)$ a.e. on $(0,1).$ Thus, we have $\Delta_a(\lambda)\equiv\Delta(\lambda),$ i.e. the spectrum of
${\cal L}_a$ as well as the one of ${\cal L}$ coincide with $\{\lambda_n\}_{n\ge 0}.$ $\hfill\Box$

\medskip
The proof of Theorem~3 is constructive and gives algorithms for solving the inverse problem. First, we
provide an algorithm for the case $\gamma\ne\pm1.$

\medskip
{\bf Algorithm 1.} Let the spectrum $\{\lambda_n\}_{n\ge 0}$ of some problem ${\cal L}(q(x),a,\gamma),$
$\gamma\ne\pm1,$ be given.

1. Construct the function $\Delta(\lambda)$ by formula (\ref{prod0}).

2. Calculate the function $w(x)\in L_2(0,1),$ inverting the Fourier transform in (\ref{FundRep}):
$$
w(x)=2\sum_{k=1}^\infty f(\pi k)\sin\pi kx,
$$
where $f(\rho)=\rho(1+\gamma^2-2\gamma\cos\rho-\Delta(\rho^2)).$

3. Find $q_a(x)\in L_2(0,1)$ by solving the non-degenerate linear system (\ref{w3}):
$$
q_a(x)=\frac{\gamma w(x)-w(1-x)}{\gamma^3-\gamma}.
$$

4. Construct the potential $q(x)$ by inverting (\ref{qa}):
\begin{equation}\label{qa->q}
q(x)=\left\{
\begin{array}{l}
\gamma q_a(x-a+1), \quad x\in(0,a),\\[2mm]
q_a(x-a), \quad x\in(a,1).\end{array}\right.
\end{equation}

\medskip
The following algorithm deals with the case $\gamma\in\{-1,1\}.$ Unlike the previous one, for definiteness it
requires specifying also an operator $K$ appeared in Theorem~2.

\medskip
{\bf Algorithm 2.} Assume that $\gamma\in\{-1,1\}$ and let the spectrum $\{\lambda_n\}_{n\ge 0}$ of a
boundary value problem ${\cal L}(q(x),a,\gamma)$ along with the operator $K$ in (\ref{dd:4-1}) with {\it
bijective} $I+\gamma K$ be given.

1. Construct the function $w(x)$ by implementing the first two steps of Algorithm~1.

2. Find the function $q_a(x)$ on the interval $(0,1/2)$ by the formula
$$
q_a\Big(\frac12-x\Big)=K\Big((I+\gamma K)^{-1}\Big(\gamma w\Big(\frac12-x\Big)\Big)\Big), \quad 0<x<\frac12.
$$

3. Find the function $q_a(x)$ on the interval $(1/2,1)$ by the formula
$$
q_a\Big(\frac12+x\Big) =\gamma w\Big(\frac12-x\Big) -\gamma q\Big(\frac12-x\Big), \quad 0<x<\frac12,
$$
and construct the potential $q(x)$ by formula (\ref{qa->q}).

\medskip
The latter algorithm also allows one to describe the set of all iso-spectral potentials $q(x),$ i.e. of those
for which the corresponding problems ${\cal L}(q(x),a,\gamma)$ (with fixed $a\in[0,1]$ and
$\gamma\in\{-1,1\})$ have one and the same spectrum $\{\lambda_n\}_{n\ge0}.$ For this purpose, on the second
step of Algorithm~2 one should use a constant operator $K,$ i.e. for which there exists a function $p(x)\in
L_2(0,1/2)$ such that
\begin{equation}\label{last}
K(f(x))=p(x)
\end{equation}
for all $f(x)\in L_2(0,1/2).$ Indeed, the following theorem holds.

\medskip
{\bf Theorem 4. }{\it If the function $p(x)$ in (\ref{last}) ranges over $L_2(0,1/2),$ then the corresponding
functions $q(x)$ constructed by Algorithm~2 form the set of all iso-spectral potentials for the given
spectrum $\{\lambda_n\}_{n\ge 0}.$}

\medskip
{\it Proof.} It is clear that, for the operator $K$ of the form (\ref{last}) and for any $p(x)\in
L_2(0,1/2),$ Algorithm~2 gives iso-spectral potentials $q(x).$ On the other hand, by virtue of Theorem~2, no
other iso-spectral potentials exist. $\hfill\Box$
\\

{\large\bf 4. Recovering from two spectra}
\\

As can be seen in the preceding section, the spectrum of the problem ${\cal L}(q(x),a,\gamma)$ possesses
sufficient (and necessary) information for unique reconstruction of the potential $q(x)$ if and only if
$\gamma\ne\pm1.$ Here, we study Inverse Problem~2, dealing with recovering $q(x)$ from both spectra
$\{\lambda_{n,\alpha}\}_{n\ge0},$ $\alpha=0,1,$ of the periodic and the antiperiodic boundary value problems
${\cal L}_\alpha(q(x),a)={\cal L}(q(x),a,(-1)^\alpha),$ $\alpha=0,1.$

Denote by $\Delta_\alpha(\lambda)$ the characteristic function of ${\cal L}_\alpha(q(x),a).$ By virtue of
Lemma~1, we have
\begin{equation} \label{FundRep_al}
\Delta_\alpha(\lambda)=2 -2(-1)^\alpha\cos\rho -\int_0^1 w_\alpha(x)\frac{\sin\rho x}\rho\,dx,\quad
w_\alpha(x)\in L_2(0,1),
\end{equation}
where the functions $w_\alpha(x)$ are determined by the formulae
\begin{equation}\label{w1}
w_\alpha(x)=\left\{
\begin{array}{l}
\displaystyle q(a+x) +q(a-x), \quad x\in(0,a),\\[3mm]
\displaystyle (-1)^\alpha q(1+a-x) +q(a+x), \quad x\in(a,1-a),\\[3mm]
\displaystyle (-1)^\alpha(q(1+a-x)+q(x-1+a)), \quad x\in(1-a,1),
\end{array}\right. \quad \alpha=0,1,
\end{equation}
which were referred as the main equations for Inverse Problem~1. Moreover, Lemma~2 gives
\begin{equation}\label{prod1}
\Delta_0(\lambda) =(\lambda-\lambda_{0,0})\prod_{n=1}^\infty\frac{\lambda_{n,0}-\lambda}{\lambda_{n,0}^0},
\quad \Delta_1(\lambda) =4\prod_{n=0}^\infty\frac{\lambda_{n,1}-\lambda}{\lambda_{n,1}^0},
\end{equation}
where $\lambda_{n,\alpha}^0,$ $n\ge0,$ $\alpha=0,1,$ are determined by (\ref{prod0-1}).

In particular, we will show that Inverse Problem~2 is uniquely solvable without any additional assumptions on
$q(x)$ if and only if $a\in\{0,1\}.$ Moreover, the following theorem along with Theorem~1 gives necessary and
sufficient conditions for its solvability in this special case.

\medskip
{\bf Theorem 5. }{\it Let $a\in\{0,1\}.$ Then for any complex sequences $\{\lambda_{n,\alpha}\}_{n\ge0},$
$\alpha=0,1,$ of the form
\begin{equation}\label{a=2}
\lambda_{2k,\alpha}=(2 k+\alpha)^2\pi^2 +\varkappa_{k,\alpha}, \;\; \{\varkappa_{k,\alpha}\}\in l_2, \;\;
k\ge0, \;\; \lambda_{2k-1,\alpha}=(2 k-\alpha)^2\pi^2, \;\; k\ge1, \;\; \alpha=0,1,
\end{equation}
there exists a unique complex-valued function $q(x)\in L_2(0,1)$ such that the given sequences are the
spectra of the corresponding boundary value problems ${\cal L}_\alpha(q(x),a),$ $\alpha=0,1.$}

\medskip
{\it Proof.} According to Remark~1, it is sufficient to consider the case $a=0.$ For $\alpha=0,1,$ construct
the entire functions $\Delta_\alpha(\lambda)$ by formulae (\ref{prod1}). According to Lemma~3, these
functions have the form (\ref{FundRep_al}) with some functions $w_\alpha(x)\in L_2(0,1)$ obeying condition
(\ref{par}), which takes the form
\begin{equation}\label{par1}
w_\alpha(x)=(-1)^\alpha w_\alpha(1-x), \quad 0<x<1, \quad \alpha=0,1.
\end{equation}
Solving the system of main equations (\ref{w1}), which, in the present case $a=0,$ takes the form
\begin{equation}\label{w4}
w_0(x)=q(1-x)+q(x), \quad w_1(x)=-q(1-x)+q(x), \quad 0<x<1,
\end{equation}
we get
\begin{equation}\label{q}
q(x)=\frac{w_0(x) +w_1(x)}2.
\end{equation}
Denote by $\tilde\Delta_\alpha(\lambda)$ the characteristic functions of the constructed problems ${\cal
L}_\alpha(q(x),0),\,\alpha=0,1.$ According to Lemma~1, they have the form
\begin{equation}\label{4.0}
\tilde\Delta_\alpha(\lambda)=2 -2(-1)^\alpha\cos\rho -\int_0^1 \tilde w_\alpha(x)\frac{\sin\rho x}\rho\,dx,
\quad \alpha=0,1,
\end{equation}
where the functions $\tilde w_\alpha(x)$ are determined by the formulae
$$
\tilde w_0(x)=q(1-x)+q(x), \quad \tilde w_1(x)=-q(1-x)+q(x), \quad 0<x<1.
$$
Comparing them with (\ref{FundRep_al}) and (\ref{w4}), we arrive at $\tilde\Delta_\alpha(\lambda)
\equiv\Delta_\alpha(\lambda),$ $\alpha=0,1,$ i.e. the spectra of the constructed problems coincide with the
given sequences, respectively.

The uniqueness follows from uniqueness of solution of the system (\ref{w4}). $\hfill\Box$

\medskip
The following algorithm gives solution of Inverse Problem~2 for $a=0.$ According to Remark~1, in the case
$a=1$ this algorithm gives $q(1-x)$ instead of $q(x).$

\medskip
{\bf Algorithm 3.} For $\alpha=0,1,$ let the spectra $\{\lambda_{n,\alpha}\}_{n\ge0}$ of the problems ${\cal
L}_\alpha(q(x),0)$ with some common potential $q(x)\in L_2(0,1)$ be given. Then it can be found by the
following steps.

1. Construct the functions $\Delta_0(\lambda)$ and $\Delta_1(\lambda)$  by formulae (\ref{prod1}).

2. Calculate the functions $w_0(x)$ and $w_1(x),$ inverting the Fourier transform in (\ref{FundRep_al}):
$$
w_\alpha(x)=2\sum_{k=1}^\infty f_\alpha(\pi k)\sin\pi kx, \quad \alpha=0,1,
$$
where $f_\alpha(\rho)=\rho(2-2(-1)^\alpha\cos\rho-\Delta_\alpha(\rho^2)).$

3. Find the function $q(x)$ by formula (\ref{q}).

\medskip
As mentioned above, we claim that 0 and 1 are unique values of $a$ in the segment $[0,1]$ for which Inverse
Problem~2 is uniquely solvable. This will actually be established by Theorem~7 below. In order to achieve the
uniqueness for $a\in(0,1),$ one should put some additional restrictions on the potential $q(x).$
Specifically, the following uniqueness theorem holds. According to Remark~1, for definiteness, we consider
the case $a\in(0,1/2].$

\medskip
{\bf Theorem 6. }{\it Let $a\in(0,1/2]$ and let there exist an operator $P:L_2(0,a)\to L_2(0,a)$ with
injective $I+P,$ such that
\begin{equation}\label{dd:4-2}
q(a-x)=P(q(a+x)), \quad \tilde q(a-x)=P(\tilde q(a+x)) \quad {\rm a.e.\;\; on} \quad (0,a).
\end{equation}
Then coincidence of the spectra: $\{\lambda_{n,\alpha}\}_{n\ge0}=\{\tilde\lambda_{n,\alpha}\}_{n\ge0},$
$\alpha=0,1,$ implies coincidence of the potentials: $q(x)=\tilde q(x)$ a.e. on $(0,1).$}

\medskip
{\it Proof.}  The coincidence of the spectra implies $w_\alpha(x)=\tilde w_\alpha(x)$ a.e. on $(0,1)$ for
$\alpha=0,1.$ Thus, by virtue of the first line in (\ref{w1}) and (\ref{dd:4-2}), we get
$$
P(q(a+x))+q(a+x)= P(\tilde q(a+x))+\tilde q(a+x) \;\; {\rm a.e. \; on} \;\: (0,a),
$$
which along with the injectivity of $I+P$ and (\ref{dd:4-2}) gives $q(x)=\tilde q(x)$ a.e. on $(0,2a).$

Further, summing up (\ref{w1}) for $x\in(a,1-a)$ and for different $\alpha$'s, we get
$$
q(a+x)=\frac{w_0(x)+w_1(x)}2 =\frac{\tilde w_0(x)+\tilde w_1(x)}2 =\tilde q(a+x), \;\; {\rm a.e. \; on} \;\:
(a,1-a),
$$
i.e. $q(x)=\tilde q(x)$ a.e. on $(2a,1),$ which finishes the proof.  $\hfill\Box$

\medskip
{\bf Remark 4.} According to Remark~1, for $a>1/2,$ condition (\ref{dd:4-2}) takes the form
\begin{equation}\label{dd:4-2-1}
q(a+x)=P(q(a-x)), \quad \tilde q(a+x)=P(\tilde q(a-x)) \quad {\rm a.e.\;\; on} \quad (0,1-a).
\end{equation}
Moreover, as in Remark~3, the operator $P$ can be applied alternatively to the left-hand sides of the
equalities in (\ref{dd:4-2}) and (\ref{dd:4-2-1}).

\medskip
The following theorem gives necessary and sufficient conditions for solvability of Inverse Problem~2 in the
case $a\in(0,1).$

\medskip
{\bf Theorem 7. }{\it Fix $a\in(0,1).$ Then for arbitrary sequences of complex numbers
$\{\lambda_{n,0}\}_{n\ge0}$ and $\{\lambda_{n,1}\}_{n\ge0}$ to be the spectra of the problems ${\cal
L}_0(q(x),a)$ and ${\cal L}_1(q(x),a),$ respectively, with a common potential $q(x)\in L_2(0,1)$ (not unique)
it is necessary and sufficient to have the form~(\ref{a=2}) and to satisfy the condition
\begin{equation}\label{growth}
\Delta_0(\lambda)+\Delta_1(\lambda)=O(\exp(\max\{1-a,a\}|\rho|)), \quad \lambda\to-\infty,
\end{equation}
where the functions $\Delta_0(\lambda)$ and $\Delta_1(\lambda)$ are determined by formulae (\ref{prod1}).}

\medskip
{\it Proof.} According to Remark~1, it is sufficient to consider the case $a\in(0,1/2].$ By necessity,
(\ref{a=2}) is already established in Theorem~1. Further, by virtue of (\ref{FundRep_al}), we get
\begin{equation} \label{D+D}
\Delta_0(\lambda) +\Delta_1(\lambda)= 4 -\int_0^1( w_0(x) +w_1(x))\frac{\sin\rho x}\rho\,dx,
\end{equation}
which along with (\ref{w1}) gives
$$
\Delta_0(\lambda) +\Delta_1(\lambda) =4 -\int_0^{1-a}( w_0(x) +w_1(x))\frac{\sin\rho x}\rho\,dx,
$$
The latter implies (\ref{growth}) and finishes the proof of the necessity.

Let some sequences $\{\lambda_{n,0}\}_{n\ge0}$ and $\{\lambda_{n,1}\}_{n\ge0},$ obeying the conditions of the
theorem be given. For $\alpha=0,1,$ construct the entire functions $\Delta_\alpha(\lambda)$ by formulae
(\ref{prod1}). According to Lemma~3, they have the form (\ref{FundRep_al}) with some functions
$w_\alpha(x)\in L_2(0,1)$ satisfying condition (\ref{par1}). Since the indicator diagrams of the functions
$\Delta_\alpha(\rho^2)$ lie on the imaginary axis of the $\rho$-plane, condition (49) implies that the type
of their sum does not exceed $\max\{1-a,a\}=1-a.$ According to the Paley--Wiener theorem, representation
(\ref{D+D}) implies
\begin{equation} \label{4.1}
w_0(x)+w_1(x)=0 \;\; {\rm a.e. \;\; on} \;\; (1-a,1).
\end{equation}
Choose any function $q(x)\in L_2(0,2a)$ that obeys the functional equation
\begin{equation} \label{4.3}
w_0(x)=q(a+1-x)+q(a-1+x), \quad x\in(1-a,1).
\end{equation}
If $a\ne1/2,$ then we consider also the linear system
\begin{equation} \label{4.4}
w_0(x)=q(1+a-x) +q(a+x), \quad w_1(x)=-q(1+a-x) +q(a+x), \quad x\in(a,1-a),
\end{equation}
which has the unique solution
\begin{equation} \label{4.5}
q(x)=\frac{w_0(x-a) +w_1(x-a)}2, \quad x\in(2a,1).
\end{equation}
Thus, we constructed some function $q(x).$ Note that (\ref{par1}) and (\ref{4.1})--(\ref{4.4}) imply
(\ref{w1}).

For $\alpha=0,1,$ consider the boundary value problems ${\cal L}_\alpha(q(x),a)$ with this $q(x).$ Denote by
$\tilde\Delta_\alpha(\lambda)$ their characteristic functions. According to Lemma~1, they have the form
(\ref{4.0}) with the functions $\tilde w_\alpha(x)$ determined by the formulae
$$
\tilde w_\alpha(x)=\left\{
\begin{array}{l}
\displaystyle q(a+x) +q(a-x), \quad x\in(0,a),\\[3mm]
\displaystyle (-1)^\alpha q(1+a-x) +q(a+x), \quad x\in(a,1-a),\\[3mm]
\displaystyle (-1)^\alpha(q(1+a-x)+q(x-1+a)), \quad x\in(1-a,1),
\end{array}\right. \quad \alpha=0,1.
$$
Comparing this with (\ref{w1}) and using (\ref{FundRep_al}) and (\ref{4.0}), we arrive at
$\tilde\Delta_\alpha(\lambda) \equiv\Delta_\alpha(\lambda)$ for $\alpha=0,1,$ i.e. the spectra of the
constructed problems coincide with the given sequences, respectively. $\hfill\Box$

\medskip
{\bf Remark 5.} According to (\ref{par1}) and (\ref{4.1}), Inverse Problem~1 for $\gamma\in\{-1,1\}$ is
equivalent to Inverse Problem~2 as soon as $a=1/2.$ Indeed, in this case, specification of the spectrum
$\{\lambda_{n,0}\}_{n\ge0}$ is equivalent to the specification of $\{\lambda_{n,1}\}_{n\ge0}.$

\medskip
By virtue of Remarks~1 and~5, it is sufficient to provide an algorithm for solving Inverse Problem~2 only for
$a\in(0,1/2).$

\medskip
{\bf Algorithm 4.} For $\alpha=0,1,$ let the spectra $\{\lambda_{n,\alpha}\}_{n\ge0}$ of the problems ${\cal
L}_\alpha(q(x),a)$ with $a\in(0,1/2]$ and $q(x)\in L_2(0,1)$ be given along with the operator $P$ in
(\ref{dd:4-2}) with {\it bijective} $I+P.$ Then the potential $q(x)$ can be found by the following steps.

1. Calculate the functions $w_0(x)$ and $w_1(x),$ implementing the first two steps of Algorithm~3.

2. Find $q(x)$ on the interval $(0,a)$ by the formula
$$
q(a-x)=P((I+P)^{-1}(w_0(x))), \quad 0<x<a.
$$

3. Find $q(x)$ on the interval $(a,2a)$ by the formula
$$
q(x+a) =w_0(x) -q(a-x), \quad 0<x<a.
$$

4. Finally, construct $q(x)$ on $(2a,1)$ by formula (\ref{4.5}).

\medskip
This algorithm allows one to describe the set of iso-bispectral potentials $q(x),$ i.e. of those for which
the corresponding problems ${\cal L}_\alpha(q(x),a),$ $\alpha=0,1,$ have one and the same pair of the spectra
$\{\lambda_{n,\alpha}\}_{n\ge0},$ $\alpha=0,1.$ For definiteness, let $a\in(0,1/2].$ Then on the second step
of Algorithm~4 one should use the constant operator
\begin{equation} \label{4.6}
P(f(x))\equiv p(x),
\end{equation}
where $p(x)\in L_2(0,a)$ is fixed.

\medskip
{\bf Theorem 8. }{\it Let $a\in(0,1/2].$ If the function $p(x)$ ranges over $L_2(0,a),$ then the
corresponding functions $q(x)$ constructed by Algorithm~4 form the set of all iso-bispectral potentials for
the given pair of spectra $\{\lambda_{n,0}\}_{n\ge 0}$ and  $\{\lambda_{n,1}\}_{n\ge 0}.$}

\medskip
{\it Proof.} It is clear that for the operator $P$ of the form (\ref{4.6}) and for any $p(x)\in L_2(0,a),$
Algorithm~4 gives iso-bispectral potentials $q(x).$ On the other hand, by virtue of Theorem~6, no other
iso-bispectral potentials exist. $\hfill\Box$
\\

{\large\bf Appendix A}
\\

Here we prove the following auxiliary assertion.

\medskip
{\bf Lemma~A1. }{\it For $\alpha\notin{\mathbb Z},$ the system
$\Lambda_\alpha:=\{\sin(2n+\alpha)x\}_{n\in{\mathbb Z}}$ is a Riesz basis in $L_2(0,\pi).$}

\medskip
For this purpose, we need the following theorem of Levin and Ljubarski\u{\i} (see \cite{LL}).

\medskip
{\bf Theorem A1. }{\it The system of functions $\{\exp(iz_k x)\}_{k\in{\mathbb N}},$ where
$\{z_k\}_{k\in{\mathbb N}}$ is the set of zeros of a sine-type function and $\displaystyle\inf_{k\ne
n}|z_k-z_n|>0,$ is a Riesz basis in $L_2(-\pi,\pi).$}

\medskip
We remind that an entire function of exponential type $S(\rho)$ is called a {\it sine-type function} if, for
some positive constants $c,$ $C$ and $K,$ the inequalities $c<|S(\rho)|\exp(-|{\rm Im}\rho|\pi) <C$ hold as
soon as $|{\rm Im}\rho|>K.$

\medskip
{\it Proof of Lemma~A1.} For real $\alpha\notin{\mathbb Z},$ this assertion was proved in \cite{Bond18}. Here
we adapt its proof to cover also the case ${\rm Im}\,\alpha\ne0.$ Firstly note, that $\Lambda_\alpha$ is
complete in $L_2(0,\pi).$ Indeed, since for any $f(x)\in L_2(0,\pi)$ we have
$$
\frac1{\sigma(\rho)} \int_0^\pi f(x)\sin\rho x\,dx = o(1), \quad \rho\to\infty, \quad \rho\in G_{\delta,1},
$$
where $G_{\delta,1}=\{\rho:|\rho-2n-\alpha|\ge \delta, n\in{\mathbb Z}\},$ $\delta>0$ and
$$
\sigma(\rho)=\sin\frac{\rho+\alpha}2\pi\cdot \sin\frac{\rho-\alpha}2\pi,
$$
then $f(x)=0$ a.e. on $(0,\pi)$ as soon as $\overline{f(x)}$ is orthogonal to all elements of
$\Lambda_\alpha.$ Thus, it remains to show the existence of positive constants $A_1$ and $A_2$ such that for
any sequence $\{c_n\}\in l_2,$ the following two-sided estimate holds:
\begin{equation} \label{A1}
A_1\sum_{n=-\infty}^\infty |c_n|^2 \le\Big\|\sum_{n=-\infty}^\infty c_n\sin(2n+\alpha)x\Big\|_{L_2(0,\pi)}^2
\le A_2\sum_{n=-\infty}^\infty |c_n|^2.
\end{equation}
Since $2n\pm\alpha,$ $n\in{\mathbb Z},$ are all zeros of the sine-type function $\sigma(\rho),$ Theorem~A1
implies that the system $\{\exp(i(2n+(-1)^\nu\alpha)x)\}_{n\in{\mathbb Z},\,\nu=0,1}$ is a Riesz basis in
$L_2(-\pi,\pi).$ Hence, there exist positive constants $B_1$ and $B_2$ such that for any $\{b_n\}\in l_2,$
the two-sided estimate
\begin{equation} \label{A2}
B_1\sum_{n=-\infty}^\infty |b_n|^2 \le\Big\|\sum_{\nu=0}^1\sum_{n=-\infty}^\infty
b_{2n+\nu}\exp(i(2n+(-1)^\nu\alpha)x)\Big\|_{L_2(-\pi,\pi)}^2 \le B_2\sum_{n=-\infty}^\infty |b_n|^2
\end{equation}
holds. In particular, with $b_{2n}=-ic_n/2$ and $b_{2n+1}=ic_{-n}/2$ for $n\in{\mathbb Z},$ estimates
(\ref{A2}) imply
\begin{equation} \label{A3}
\frac{B_1}2\sum_{n=-\infty}^\infty |c_n|^2 \le\Big\|\sum_{n=-\infty}^\infty
c_n\sin(2n+\alpha)x\Big\|_{L_2(-\pi,\pi)}^2 \le \frac{B_2}2\sum_{n=-\infty}^\infty |c_n|^2.
\end{equation}
Since $\|f\|_{L_2(-\pi,\pi)}=\sqrt2\|f\|_{L_2(0,\pi)}$ for any odd (as well as for any even) function
$f(x)\in L_2(-\pi,\pi),$ estimates (\ref{A3}) give (\ref{A1}) with $A_j=B_j/4,$ $j=1,2,$ which finishes the
proof. $\hfill\Box$
\\

{\bf Acknowledgement.} The authors are grateful to Professor Natalia Bondarenko for sharing the idea of
proving Lemma~A1.

\medskip
{\bf Funding.} The first author was supported by Grant 20-31-70005 of the Russian Foundation for Basic
Research.

\medskip
{\bf Availability of data and material.} The manuscript has no associated data.

\end{document}